\documentclass[final]{siamltex}
\usepackage{amsmath}
\usepackage{amsfonts}
\usepackage{epsfig}

\title{Spectral Methods for Partial Differential \\
Equations in Irregular Domains: \\
The Spectral Smoothed Boundary Method \thanks{This work was supported by 
grants BFM2003-02832 (Ministerio de Ciencia y Tecnolog\'{\i}a, Spain), PAC-02-002 (Consejer\'{\i}a de Ciencia y Tecnolog\'{\i}a 
de la Junta de Comunidades de Castilla-La Mancha -CCyT-JCCM-, Spain)  .}}

\author{Alfonso Bueno-Orovio\thanks{Departamento de Matem\'aticas,
Universidad de Castilla-La Mancha,
E.T.S.I.  Industriales, Av. Camilo Jos\'e Cela, 3, Ciudad Real,
E-13071, Spain (alfonso.bueno@uclm.es). Supported by 
 CCyT-JCCM under PhD grant \mbox{03-056}. }
\and V\'{i}ctor M.
P\'erez-Garc\'{i}a \thanks{Departamento de Matem\'aticas,
Universidad de Castilla-La Mancha,
E.T.S.I.  Industriales, Av. Camilo Jos\'e Cela, 3, Ciudad Real,
E-13071, Spain (victor.perezgarcia@uclm.es).} \and Flavio H.
Fenton\thanks{Department of Physics, Hofstra University and Beth
Israel Medical Center, NY 1003, USA (Flavio.Fenton@hofstra.edu).}}

\begin{document}%SourceDoc 

\maketitle \vspace{-1.35in} %\slugger{sisc}{xxxx}{xx}{x}{xxxx-xxxx}
\vspace{1.1in}

\setcounter{page}{1}

\begin{abstract}
In this paper, we propose a numerical method to approximate the solution of partial differential equations
in irregular domains with no-flux boundary conditions by means of spectral methods. The main features
of this method are its capability to deal with domains of arbitrary
shape and its easy implementation via Fast Fourier Transform routines. We discuss several examples 
of practical interest and test the results against exact solutions and standard numerical methods.
\end{abstract}



\pagestyle{myheadings} \thispagestyle{plain} \markboth{A.
BUENO-OROVIO, V. M. P\'EREZ-GARC\'IA AND F. H.
FENTON}{SPECTRAL METHODS IN IRREGULAR DOMAINS}

\section{Introduction}

Spectral methods \cite{Fornberg,Trefethen,Sanzserna} are among the most extensively used methods for the
discretization of spatial variables in partial differential equations
and have been shown to provide
very accurate approximations of sufficiently smooth solutions. Because of their high-order accuracy, the use of
spectral methods has become widespread over the years in various fields, including fluid dynamics, quantum
mechanics, heat conduction and weather prediction \cite{Canuto,Gottlieb,Roger,solitones,BunYuGuo}. However, these
methods have some limitations which have prevented them from being extended to many problems where
finite-difference and finite-element methods continue to be used predominantly. One limitation is that the discretization of
partial differential equations by spectral methods leads to the solution of large systems of linear or nonlinear
equations involving \emph{full} matrices. Finite-difference and finite-element methods, on the other hand, lead
to systems involving sparse matrices that can be handled by appropriate methods to reduce the complexity of the
calculations substantially. Another drawback of spectral methods is that the geometry of the problem domain must
be simple enough to allow the use of an appropriate orthonormal basis to expand the full set of possible
solutions to the problem. This inability to handle irregularly shaped domains is one reason why these methods
have had limited use in many engineering problems, where finite-element methods are preferred because of their
flexibility to describe complex geometries despite the computational costs associated with constructing an
appropriate solution grid. Although there have been attempts to use spectral methods in irregular domains
\cite{Orszag,Korczak}, these approaches usually involve either incorporating finite-element preconditioning or 
the use of so-called spectral elements which are similar to finite elements. We are
not aware of any previous study where purely spectral methods, particularly those involving Fast Fourier Transforms (FFTs), have been used to
obtain solutions in complex irregular geometries.

In this paper we present an accurate and easy-to-use method for 
approximating the solution of partial
differential equations in irregular domains with no-flux boundary conditions using spectral methods. The idea is
based on what in dendritic solidification is known as the phase-field method \cite{KarRap98}. This method is used to
locate and track the interface between the solid and liquid states and has been applied to a wide
variety of problems including viscous fingering \cite{Foletal99a,Foletal99b}, crack propagation \cite{Karetal01}
and the tumbling of vesicles \cite{BibMis03}. For a comprehensive review see \cite{Casademunt}.

In what follows we use the idea behind phase-field
methods to illustrate how the solution of several partial
differential equations can be obtained in various irregular and
complex domains using spectral methods. Throughout the manuscript,
for simplicity, we will refer to the combination of the
phase-field and spectral methods as the spectral smoothed
boundary (SSB) method. Our approach consists of two steps. First,
the idea of the phase-field method is formalized and its
convergence analyzed for the case of homogeneous Neumann boundary
conditions. Then we discuss how the new formulation is useful for
the direct use of spectral methods, specifically those based on
trigonometric polynomials. This formulation makes the problem
suitable for efficient solution using FFTs \cite{FFTW}. Since it
is our intention that the resulting methodology be used in a
variety of problems in engineering and applied science, we have
concentrated on the important underlying concepts, reserving some
of the more formal questions related to these methods for a
subsequent analysis.

\section{The phase-field (smoothed boundary) method}
\label{ppi} In this work we focus on applying the
phase-field method to
partial differential equations of the form
\begin{subequations}
\begin{equation}\label{equa}
\nabla (\boldsymbol{D}^{(j)}\nabla u_j) + f(u_1,...,u_N,t) =
\partial_t u_j
\end{equation}
for $N$ unknown real functions $u_j$ defined on an
irregular domain $\Omega \subset \mathbb{R}^n$, where
$n=1,2,3$ is the spatial dimensionality of the problem,
together with appropriate initial conditions
$u_j(x,0) = u_{j0}(x)$ and subject to Neumann boundary conditions
 \begin{equation} \label{boundary}
 (\boldsymbol {\Vec{n}} \cdot \boldsymbol{D}^{(j)}\nabla u_j) = 0
 \end{equation}
 \end{subequations}
on $\partial \Omega$, where $\boldsymbol{D}^{(j)}(x)$ is
a family of $n\times n$ matrices that may depend on the spatial variables.
 Equations (\ref{equa}) and (\ref{boundary}) include many
reaction-diffusion models, such as those describing population dynamics
or cardiac electrical activity. In here we will restrict the analysis to equations of the form \eqref{equa} although we believe that the idea behind the method can be extended to 
many other problems involving complex boundaries and different types of partial differential equations.

Instead of discretizing Eq.
(\ref{equa}) the smoothed boundary method relies on considering the auxiliary problem
\begin{equation}\label{equa2}
\nabla (\phi^{(\xi)} \boldsymbol{D}^{(j)}\nabla u^{(\xi)}_j)
+ \phi^{(\xi)}f(u^{(\xi)}_1,...,u^{(\xi)}_N,t) =
\partial_t (\phi^{(\xi)}u^{(\xi)}_j),
\end{equation}
for the unknown functions $u^{(\xi)}_j$ on an enlarged domain 
 $\Omega'$ satisfying the following conditions: (i) 
$\Omega \subset \Omega'$ and (ii) $\partial \Omega \cup \partial \Omega' = \emptyset$.  
The function $\phi^{(\xi)}$ is continuous  in  $\Omega'$ and has the value one inside
$\Omega$ and
smoothly decays to zero outside $\Omega$, with  $\xi$ identifying
the width of the decay. That is, if $\chi_{\Omega}$
is the characteristic function of the set $\Omega$ defined as
\begin{equation}\label{chieq}
\chi_{\Omega}=\left\{\begin{array}{cl}
    1, & x \in \Omega \\
    0, & x \in \Omega'-\Omega
       \end{array}\right.
\end{equation}
then $\phi^{(\xi)}: \Omega' \rightarrow \mathbb{R}$ is a
regularized approximation to $\chi_{\Omega}$ such that
$\lim_{\xi\rightarrow 0} \phi_{\xi} = \chi_{\Omega}$.

\begin{figure}
\begin{center}
\epsfig{file=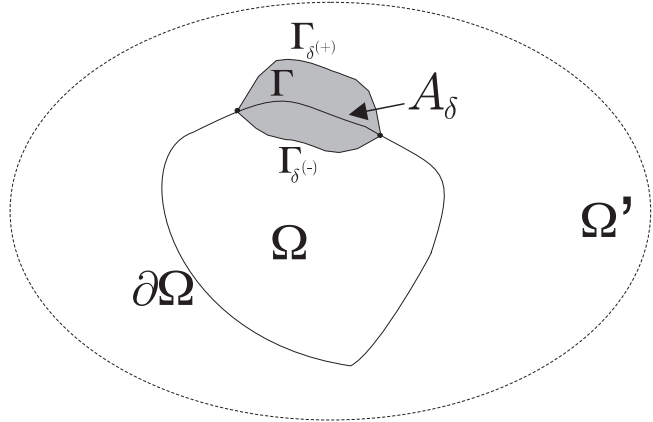,width=0.66\textwidth}
\caption{Illustration of an irregular domain $\Omega$ 
with an example of the connected curve
$\Gamma$ and the domain $A_\delta$ used in the proof of convergence of the method.
\label{prima}}
\end{center}
\end{figure}

The key idea of the Smoothed boundary method (SBM) is that  when $\xi \rightarrow 0$ the solutions $u_j^{(\xi)}$ of Eqs. (\ref{equa2})
on any domain $\Omega'$ with arbitrary boundary conditions on $\partial \Omega'$ satisfy
 $u^{(\xi)}_j \rightarrow u_j$, that means, they tend to the solution of Eqs. (\ref{equa}), automatically incorporating the boundary conditions (\ref{boundary}). To see why, let us first realize that inside $\Omega$ the statement is inmediate since 
 $\phi^{(\xi)}\rightarrow 1$ in $\Omega$ as $\xi
\rightarrow 0$ and Eq. \eqref{equa2} becomes Eq. \eqref{equa}.
At the boundary we consider for simplicity the
situation with $n=2$ (the extension to $n=3$ is immediate).
Assuming smoothness of $\partial \Omega$ (which in
this case will be a curve) and choosing any connected curve $\Gamma \subset
\partial \Omega$, we define two families of differentiable curves $\Gamma_{\delta^{(+)}} \subset \Omega'/\Omega$ and
 $\Gamma_{\delta^{(-)}} \subset \Omega$  whose ends coincide with those of $\Gamma$ 
 and which tend uniformly to $\Gamma$ following the parameter $\delta$.
 The curves $\Gamma_{\delta^{(+)}}$ and
 $\Gamma_{\delta^{(-)}}$ are then the boundaries of a region 
 $A_{\delta}$ whose boundary
$\partial A_{\delta}=\Gamma_{\delta^{(+)}} \cup
\Gamma_{\delta^{(-)}}$ (see Fig. \ref{prima}).

We now integrate Eq. (\ref{equa2}) over $A_{\delta}$:
\begin{equation}\label{integr}
\int \int_{A_{\delta}} \left[ \nabla (\phi^{(\xi)}
\boldsymbol{D}^{(j)}\nabla u^{(\xi)}_j) +
\phi^{(\xi)}f(u,t) \right] dx = \int \int_{A_{\delta}}
\partial_t (\phi^{(\xi)}u^{(\xi)}_j) dx.
\end{equation}
Using the Gauss (or Green) theorem for the first term of Eq.
(\ref{integr}) we obtain
\begin{equation}\label{intro}
 \oint_{\partial A_{\delta}} \boldsymbol{n} \cdot \phi^{(\xi)}
\boldsymbol{D}^{(j)}\nabla u^{(\xi)}_j dx + \int \int_{A_{\delta}}\phi^{(\xi)}f(u,t) dx = \int \int_{A_{\delta}}
\partial_t (\phi^{(\xi)}u^{(\xi)}_j) dx,
\end{equation}
where $\oint$ denotes a line integral over $\partial A_{\delta}$.

Now we take the limit $\xi \rightarrow 0$ in (\ref{intro})
 to obtain
\begin{eqnarray}
\lim_{\xi \rightarrow 0} \oint_{\partial A_{\delta}} \boldsymbol{n} \cdot \phi^{(\xi)}
\boldsymbol{D}^{(j)}\nabla u^{(\xi)}_j dx & = & - \lim _{\xi \rightarrow 0} \int
\int_{A_{\delta}}\phi^{(\xi)}f(u,t) dx \nonumber
\\ & & + \lim_{\xi \rightarrow 0} \int  \int_{A_{\delta}}
\partial_t (\phi^{(\xi)}u^{(\xi)}_j) dx \nonumber \\
& & = m(A_{\delta}) \left[ -\phi^{(\xi)}f(u,t) +
\partial_t (\phi^{(\xi)}u^{(\xi)}_j) \right]_{x=\zeta},
\label{finita}
\end{eqnarray}
where the last equality comes from the mean value theorem for integrals and $m(A_{\delta})$ is the measure of
the set $A_{\delta}$. Here we assume that the solutions to Eq. \eqref{equa2} and its time derivatives are
bounded so that the right-hand side of Eq. (\ref{finita}) is finite. On the left-hand side we decompose
$\oint_{\partial A_{\delta}}$ as $\int_{\Gamma_{\delta^{(+)}}} + \int_{\Gamma_{\delta^{(-)}}}$. It is evident that $\lim_{\xi \rightarrow 0} \int_{\Gamma_{\delta^{(+)}}} \boldsymbol{n} \cdot
\phi^{(\xi)} \boldsymbol{D}^{(j)}\nabla u^{(\xi)}_j dx = 0$, since in this limit $\phi^{(\xi)}=0$ over all $\Gamma_{\delta^{(+)}}$ and $\phi^{(\xi)}=1$
on $\Gamma_{\delta^{(-)}}$. 
As we were interested in proving that the boundary conditions are satisfied, we now make the width of the integration region $A_\delta$ tend to zero. Since $\Gamma_{\delta(-)} \rightarrow  \Gamma$ as $\delta \rightarrow 0$, and $\lim_{\delta \rightarrow 0} m(A_{\delta}) = 0$, we obtain
\begin{eqnarray}\label{final}
\lim_{\delta\rightarrow 0} \lim_{\xi\rightarrow 0} 
\oint_{\partial A_{\delta}} \boldsymbol{n} \cdot \phi^{(\xi)}
\boldsymbol{D}^{(j)}\nabla u^{(\xi)}_j dx & = & \lim_{\delta\rightarrow 0} \int_{\Gamma_{\delta(-)}} \boldsymbol{n} \cdot \boldsymbol{D}^{(j)}\nabla u^{(\xi)}_j dx \nonumber \\
& = & \int_{\Gamma}
\boldsymbol{n} \cdot \boldsymbol{D}^{(j)}\nabla u^{(\xi)}_j dx = 0.
\end{eqnarray}
Since Eq. (\ref{final}) is true for any boundary segment $\Gamma$, we obtain the final
result that in the limit $\xi \rightarrow 0$ Eq. \eqref{equa2}
satisfies $\boldsymbol{n} \cdot \boldsymbol{D}^{(j)}\nabla u^{(\xi)}_j =
0$ for $j=1,...,N$, i.e., the boundary conditions.

From what we have shown it is clear that $\Omega'$ can be any
closed region containing $\Omega$. The idea of the smoothed
boundary method is then to consider Eq. \eqref{equa2} for a small
but finite $\xi$ and to discretize this problem instead of Eqs.
(\ref{equa}) and (\ref{boundary}). The main advantage is that one
can search for the approximation $u^{(\xi)}$ on any enlarged domain
$\Omega'$ such that $\Omega \subset \Omega'$ (for instance, a
rectangular region). The enlarged discrete problem can then be solved with
any proper boundary conditions on $\partial \Omega'$, since the
fulfilment of the boundary conditions for $u$ on $\partial \Omega$
is guaranteed in the limit $\xi \rightarrow 0$. In our case, we
will use the basis of trigonometric polynomials $e^{ikx}$
to approximate the solutions; thus, we will seek an extension of
the solution $u$ of Eqs. \eqref{equa} and  \eqref{boundary}  that is
periodic on the enlarged region $\Omega'$.

\section{The Spectral Smoothed Boundary Method}

We want to discretize Eq. (\ref{equa2}) on an enlarged region. As discussed earlier, 
we choose $\Omega'$ to be a rectangular region containing $\Omega$ and we will expand $u^{(\xi)}$
in the basis of Cartesian products of trigonometric polynomials $e^{ik_xx}e^{ik_yy}$.

Let us rewrite Eq. (\ref{equa2}) without subscripts and superscripts as
\begin{equation} \label{equa3}
\nabla \phi \cdot \boldsymbol{D}\nabla u
+ \phi \nabla (\boldsymbol{D}\nabla u) + \phi f(u,t) =
\partial_t (\phi u).
\end{equation}
Note that since $\phi$ is located inside of the time derivative of
the right term of Eq. (\ref{equa2}), it is possible for the
integration domain itself to evolve in time, and thus this method
could be used to solve moving boundary problems once a coupling
equation is added for the movement of $\phi$. However, in this manuscript  we
will only deal with stationary integration domains; thus,  $\partial_t \phi = 0$  and the right side of
Eq. (\ref{equa3}) can be simplified as $\phi \partial_t u$. Dividing  Eq. (\ref{equa3}) by $\phi$, we get
\begin{equation} \label{equa4}
\nabla \log \phi \cdot \boldsymbol{D}\nabla u
+ \nabla (\boldsymbol{D}\nabla u) + f(u,t) =
\partial_t u,
\end{equation}
which is the equation of the smoothed boundary method that we will
use to perform numerical simulations.

\begin{figure}
\begin{center}
\epsfig{file=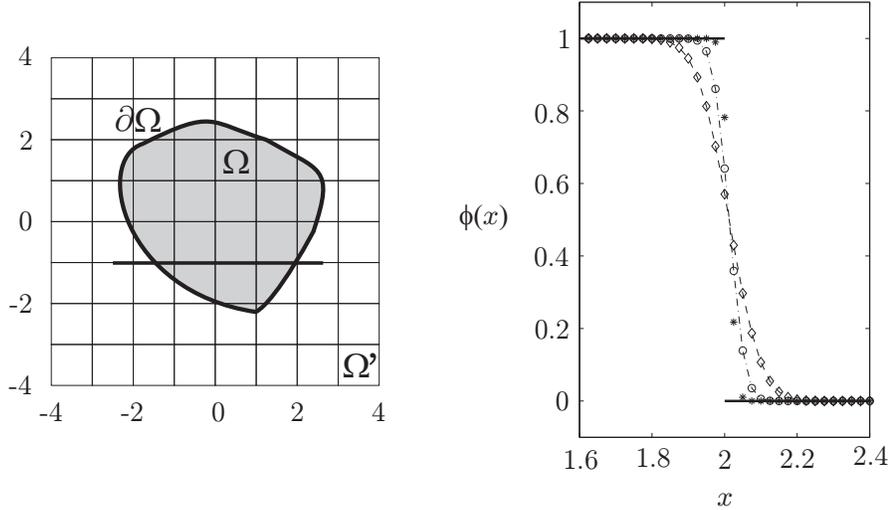,width=0.90\textwidth} 
\caption{Left: 
Example of an irregular domain $\Omega$ defined on a Cartesian 
grid and an enlarged domain $\Omega'$.
Right: Smoothing of the irregular boundary in a one-dimensional section of 
the domain. The solid line shows a small section of the 
characteristic function  $\chi_{\Omega}$ 
(with value 0 or 1) corresponding to part of the thicker line shown in the left part of the
figure.  
Phase-field functions $\phi^{(\xi)}$ obtained from  $\chi_{\Omega}$
for $\xi = 0.10$, $0.05$ and $0.025$ are labeled by
diamonds, circles and stars respectively.}
\label{fig:phi}
\end{center}
\end{figure}

To implement numerically any solution method for Eq. \eqref{equa4}, we need to 
make a specific choice for $\phi^{(\xi)}$. In practice, 
any method that produces a smooth characteristic function can be used.
In the context of phase-field
methods, the standard procedure for obtaining the values of
$\phi^{(\xi)}$ (which is called the ``phase-field") is to
integrate an auxiliary diffusion equation of the form
$\partial_t \phi =  \xi^2 \nabla^2 \phi +(2 \phi-1)/2 - (2
\phi-1)^3/2, \label{evalphi}$
with initial conditions $\phi^{(\xi)}(t=0)$ = $\chi_{\Omega}$,
until a steady state is reached \cite{Casademunt,Fenton}.
Alternatively, since we only seek a smoothed boundary we choose to obtain
$\phi^{(\xi)}$ from $\chi_{\Omega}$ using a convolution
of the form
\begin{equation}\label{convolution}
\phi^{(\xi)} = \chi_{\Omega} \ast G^{({\xi})},
\end{equation}
where $G^{(\xi )}$ is any family of functions such
that $\lim_{\xi \rightarrow 0} G^{(\xi)}(x) = \delta(x)$, where
$\delta$ is the Dirac delta function. In particular
Gaussian functions of the form
\begin{equation} \label{gaussian}
G^{(\xi)}(x) = \prod_{k=1}^n \exp(-x_k^2/\xi^2)
\end{equation}
can be chosen. In this paper, all the functions $\phi^{(\xi)}$ have been obtained
using this $n$-dimensional discrete convolution of
$\chi_{\Omega}$ with a Gaussian function of the form given by
Eq. (\ref{gaussian}). An example of the creation of
$\phi^{(\xi)}$ is
shown in Fig. \ref{fig:phi}, where it can be seen that
the width of the interface in which $\phi^{(\xi)}$ changes from zero to one 
depends on the value used for $\xi$ (in fact it is of order $\xi$).

To avoid computational difficulties for very small values of $\phi$ we approximate
 $\log \phi \approx \log (\phi + \epsilon)$, where $\epsilon$ is the machine precision.
Numerically, $\phi$ and $(\phi + \epsilon)$ are equal up to roundoff errors, but this correction bounds the
value of $\log \phi$ as $\phi \rightarrow 0$. Choosing $\phi$ to be a periodic function to avoid Gibbs phenomenom when computing its derivatives forces us to choose a computational domain in which this function becomes small enough near $\partial \Omega'$. In practice, this restriction requires us to leave a reasonable margin between the boundaries of the
physical and the enlarged domains. We found that a
margin of value $M = 10\xi$ is sufficient for all the simulations
 to be stable.

All the spatial derivatives in Eq. (\ref{equa4}) are computed in Cartesian coordinates with spectral accuracy in Fourier space. If $g(x)$ is a periodic and sufficiently smooth function, then its $n^{th}$ derivative is given by
\begin{equation}
\frac{\partial^n g}{\partial x^n} =
\mathcal{F}^{-1}\left\{ (ik_x)^n \mathcal{F} \{g\} \right\},
\end{equation}
where $\mathcal{F}$ and $\mathcal{F}^{-1}$ 
denote the direct and inverse Fourier Transform respectively, $k_x$ are the wave numbers
associated with each Fourier mode, and $i$ is the imaginary
unit.  As mentioned previously, the use of this representation for $u$ implicitly assumes
 periodic boundary conditions on $\partial \Omega'$.
It is significant that only Fourier Transforms are used for
these calculations instead of differentiation matrices,
thereby avoiding the generation and storage of these matrices
and yielding more efficient codes and shorter execution
times, especially when Fast Fourier Transforms routines are used.

In this paper we are not concerned with designing the most efficient 
SSBM, but only to prove that such a method can be used to integrate PDEs in irregular domains. Thus,
for time integration we use a simple second-order explicit
method. In the
particular case where all the coefficients of the diffusion tensor
$\boldsymbol{D}$ are constants, we can write
Eq. (\ref{equa4}) in the form
\begin{equation}
\mathcal{L}u + \mathcal{N}(u,t) = \partial_t u,
\end{equation}
where $\mathcal{L}u = \nabla (\boldsymbol{D}\nabla u)$ is the
linear term and $\mathcal{N}(u,t) = \nabla \log \phi \cdot
\boldsymbol{D}\nabla u + f(u,t)$ is the nonlinear part of the
equation. Then a second-order in time operator splitting scheme of the
form
\begin{equation}
U(t + \Delta t)  = e^{\mathcal{L} \Delta t/2}
e^{\mathcal{N} \Delta t} e^{\mathcal{L} \Delta t/2} U(t)
\end{equation}
can be used to solve the equation in time \cite{Strang}. For the examples to be presented later, we solve the nonlinear term by a second-order (half-step) explicit
method and integrate the linear part exactly in Fourier space by
exponential differentiation, which reduces the stiffness of the
problem considerably and allows the use of larger time steps. The operator splitting scheme can then be written as
\begin{subequations}
\begin{eqnarray}
U^{\ast} & = & \mathcal{F}^{-1} \{e^{\mathcal{L}\Delta t/2} \mathcal{F} \{ {U}^k \} \} \\
U^{\ast \ast} & = & U^{\ast} + \mathcal{N}(U^{\ast} +  \mathcal{N}(U^{\ast},t_k + \Delta t/2) \cdot \Delta t/2, t_k + \Delta t) \cdot \Delta t \\
U^{k+1} & = & \mathcal{F}^{-1} \{e^{\mathcal{L}\Delta t/2} \mathcal{F} \{ {U}^{\ast \ast} \} \}.
\end{eqnarray}
\end{subequations}

\section{Examples of the methodology}
 \label{applications}

\subsection{The heat equation}

\begin{figure}
\begin{center}
\epsfig{file=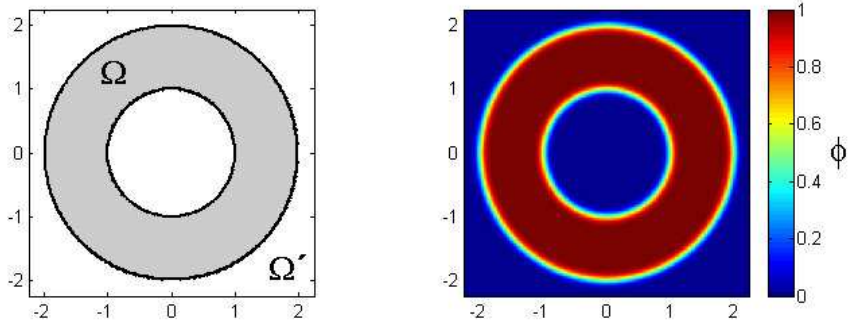,width=0.90\textwidth}
\caption{Left: The rectangular domain $\Omega'$ in which 
Eq. (\ref{eqheat}) is solved using the SSB method, with the irregular 
domain $\Omega$ (an annulus) shown in gray. 
Right:
Smoothed boundary function $\phi^{(\xi)}$ given by Eq. \eqref{convolution} for $\xi = 0.10$.}
\label{fig:phi2d}
\end{center}
\end{figure}

As a first example, we will consider a simple linear heat equation. This first case will allow us to make a quantitative study of the errors of the SSB method. Specifically, we are interested in solving the following
heat equation with sources:
\begin{equation} \label{eqheat}
\partial_t u = D \Delta u - r \cos(2\theta)
\end{equation}
in the annulus $\Omega$ defined by $1 \leq r \leq 2$ with homogeneous Neumann boundary conditions on $\partial
\Omega$, $\partial_r u |_{r = 1} = \partial_r u |_{r = 2} = 0$, and initial data $u(r,\theta,0) = 0$. The
diffusion coefficient is taken to be constant with value $D = 1$. Figure \ref{fig:phi2d} shows one example of
the generation of the smoothed boundary for this domain.

\begin{figure}
\begin{center}
\epsfig{file=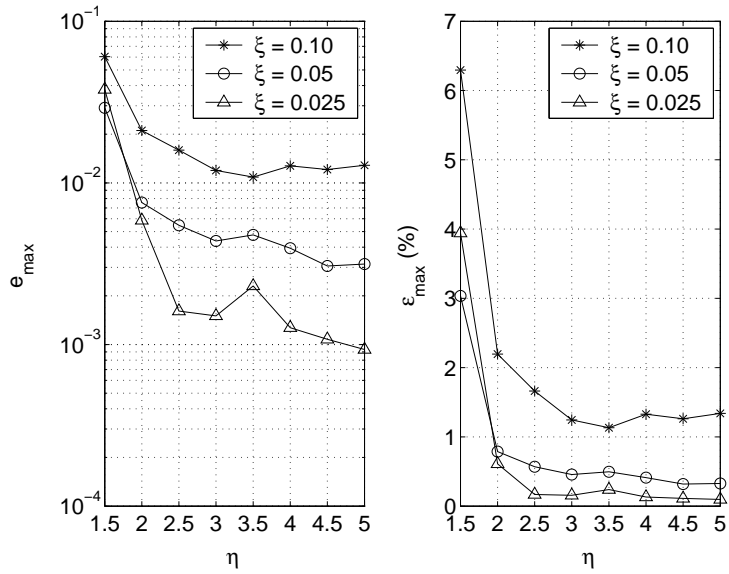,width=0.75\textwidth} 
\caption{
Maximum absolute (left) and relative
(right) errors of the numerical solution of 
Eq. (\ref{eqheat}) at time $t = 6$ in the annulus 
compared to the steady solution 
%($e_{max}/\max(u_{st})$) are shown 
as a function of the parameter 
 $\eta=\xi/\Delta x$. }
\label{fig:errorcircle}
\end{center}
\end{figure}

Equation (\ref{eqheat}) in this geometry has an explicit steady state-solution of the form
\begin{equation} \label{eqsteady}
u_{st}(r,\theta) = \left( \frac{1}{5}r^3 - \frac{31}{50}r^2
- \frac{8}{25} \frac{1}{r^2} \right) \cos(2\theta),
\end{equation}
which can be compared with the numerical steady solution of Eq. \eqref{eqheat} (in practice, we stop simulations at $t = 6$ since by this time the numerical solutions have approximately reached the steady state) and obtain error estimates. Figure \ref{fig:errorcircle}
shows the maximum absolute error $E = \| u - U^{(\xi)}\|_{\infty}$ and relative error  $e = \| u - U^{(\xi)} \|_{\infty}/\|u\|_{\infty}$  of
several simulations for different values of $\xi$ and
grid resolutions, where the relative error is defined with respect to the maximum value of the analytical solution. Note that in general these maximum errors decrease as the thickness
of the interface is reduced ($\xi \rightarrow 0$), and
that in most of the simulations the relative error is less than $1\%$.

\begin{figure}
\begin{center}
\epsfig{file=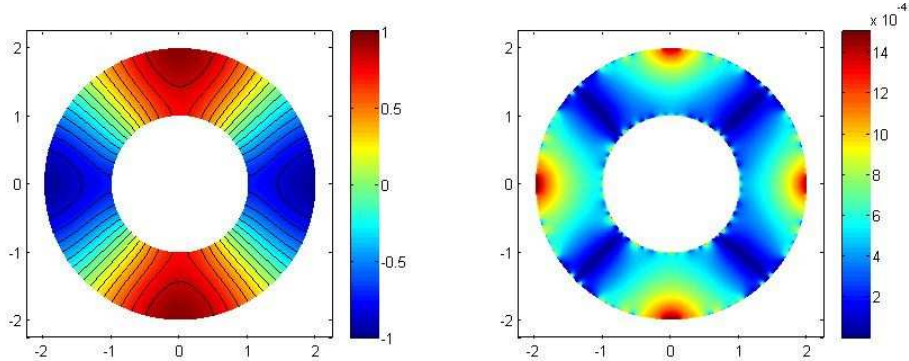,width=0.95\textwidth}
\caption{(Left) Solution of Eq. (\ref{eqheat}) 
obtained with the SSB method at time $t = 6$ 
in the annulus \mbox{$1 \le r \le 2$}. 
Grid size is $540 \times 540$ with $\xi =
0.025$. (Right) Spatial distribution of the absolute error \mbox{$E = \| u - U^{(\xi)}\|_{\infty}$} over the annulus for the simulation 
solution 
compared to the analytical steady solution. }
\label{fig:heatcircle}
\end{center}
\end{figure}

Although $\phi^{(\xi)}$ is a continuous function, it is necessary to 
have a grid fine enough to resolve properly the boundary layers in which it quickly
changes from zero to one.  For this reason, errors are
represented in Fig. \ref{fig:errorcircle} not as a function of the
number of grid points but as a function of the parameter $\eta =
\xi/\Delta x$, which gives an idea of the number of points that
lie in the interface. As mentioned before, it is also necessary to use
a margin between the integration
domain $\Omega$ and the computational domain $\Omega'$ for $\phi$
to become sufficiently small on $\partial\Omega'$. In this particular case,
\begin{equation}
\eta = \frac{\xi}{\Delta x} = \xi \frac{N}{2(R+M)},
\end{equation}
where $R$ is the outer radius of the annulus, $N$
is the grid resolution, and $M = 10\xi$ is the margin used.
This expression implies that for solving Eq. (\ref{eqheat})
in the range $\eta \in [1.5,5]$, the grid
resolution varies from 90 to 300 points if $\xi = 0.10$,
from 150 to 500 points if $\xi = 0.05$, and
equivalently from 270 to 900 points when $\xi = 0.025$.
The last important point concerning Fig.
\ref{fig:errorcircle} is that, once the interface is
properly solved ($\eta \sim 3-4$), the error converges to
an approximately constant value that depends only on
$\xi$, so there is no reason for using excessively
clustered grids. To ensure that this error is produced only
by the spatial discretization and is not due to the order of
the method chosen to perform the time integration, we
have also run the simulations with a first-order explicit
(Euler) time-integration method and obtained errors of
the same order of magnitude. Figure \ref{fig:heatcircle}
shows the solution to Eq. (\ref{eqheat}) obtained at time
$t = 6$ with the SSB method. Contour lines are also
included to illustrate that the no-flux boundary conditions
at $r = 1$ and $r = 2$ are satisfied.

\begin{figure}
\begin{center}
\epsfig{file=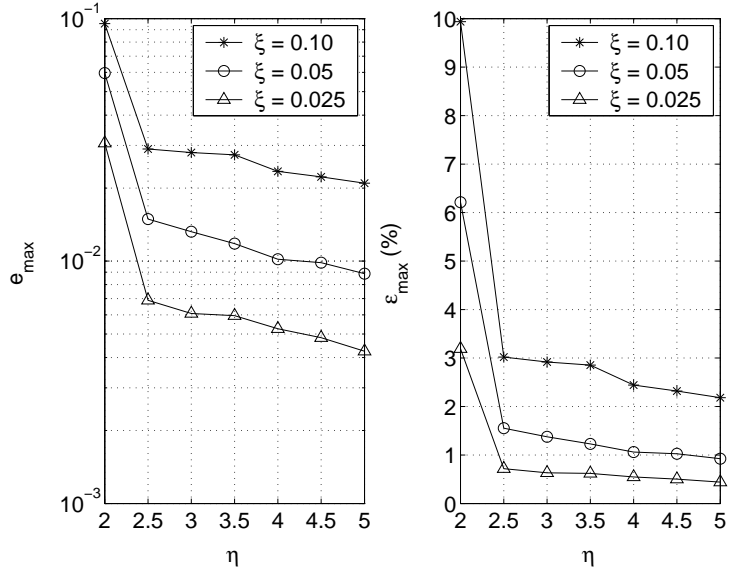,width=0.75\textwidth} 
\caption{Maximum absolute (left) and relative (right) 
errors of the numerical solution of Eq.
(\ref{eqheat}) in the quarter-annulus
$1 \leq r \leq 2$, $0 \leq \theta \leq \pi/2$ at time $t = 6$ compared to the steady
solution of the problem
as a function of the parameter 
 $\eta=\xi/\Delta x$.
} \label{fig:errorquarter}
\end{center}
\end{figure}

\begin{figure}
\begin{center}
\epsfig{file=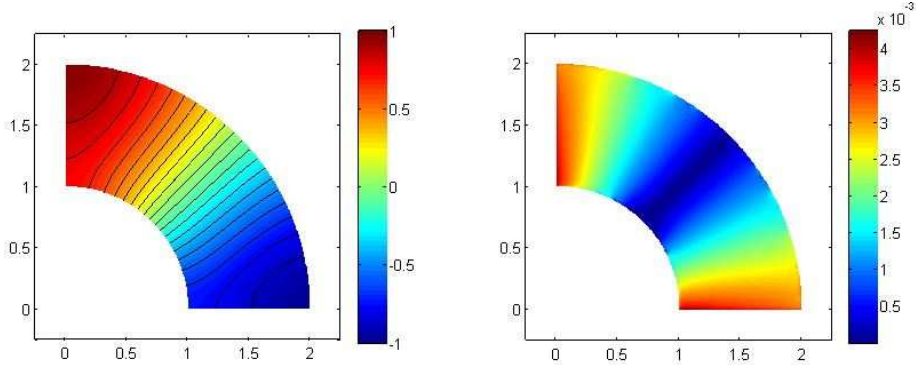,width=0.95\textwidth} 
\caption{(Left) Solution to Eq. (\ref{eqheat}) in the quarter annulus 
$1 \leq r \leq 2$, $0 \leq \theta \leq \pi/2$ at $t = 6$ using the SSB method.
Grid resolution is $500 \times 500$ and $\xi = 0.025$.
(Right)
Spatial distribution of absolute error $E = \| u - U^{(\xi)}\|_{\infty}$ over the quarter-annulus 
for the simulation solution compared to the 
analytical steady solution.
} \label{fig:heatquarter}
\end{center}
\end{figure}

The analytical solution of Eq. (\ref{eqheat}) also
satisfies homogeneous Neumann boundary conditions on the
quarter-annulus delimited by $1 \leq r \leq 2$,
$0 \leq \theta \leq \pi/2$ (see Fig. \ref{fig:heatcircle}),
which allows us to use this related geometry to show
how the SSB method performs when sharp corners are present in a given geometry.
Maximum absolute and relative errors of the simulations
for the quarter-annulus are shown in Fig.
\ref{fig:errorquarter}. Both geometries show similar good
convergence properties. However, errors are slightly larger but of the same order of magnitude than
for the full annulus due to the presence of the sharp
corners, which become slightly blunted. This can be seen
in Fig. \ref{fig:heatquarter}, which shows the error
distribution $E = \| u - U^{(\xi)}\|_{\infty}$ over the domain $\Omega$, along with the
corresponding solution.

\begin{figure}
\begin{center}
\epsfig{file=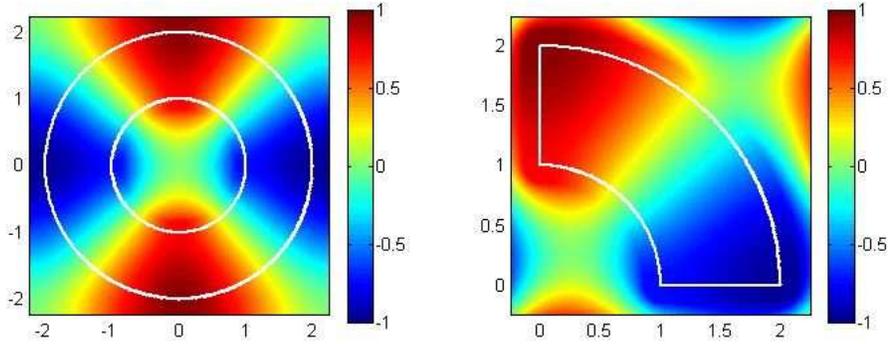,width=0.95\textwidth} 
\caption{Solution to Eq. (\ref{eqheat}) using the SSB method over the entire
domain of integration $\Omega'$ for the annulus shown in Fig.\ref{fig:heatcircle}
(left) and the quarter-annulus shown in Fig.\ref{fig:heatquarter} (right). 
Note that the periodicity along $\partial \Omega'$ imposed 
by the FFT does not alter the solution in $\Omega$ or at the boundary
$\partial\Omega$ (shown in white) where the zero-flux boundary conditions are satisfied.
} \label{fig:periodic}
\end{center}
\end{figure}

In Figs. \ref{fig:heatcircle} and \ref{fig:heatquarter} we have shown the solutions of Eq. (\ref{eqheat}) within the irregular
geometries $\Omega$. However,  the solutions $U^{(\xi)}$ are calculated over the entire domain $\Omega'$. Figure 
\ref{fig:periodic} shows the solutions over $\Omega'$ in the full and quarter-annulus examples. While no-flux 
boundary conditions are implemented along $\partial \Omega$, the overall solution has
periodic boundary conditions.  Note that as the solution $U^{(\xi)}$ is not discontinuous on $\partial \Omega$, our solutions never present Gibbs phenomena due to the irregular boundaries.

\begin{figure}
\begin{center}
\epsfig{file=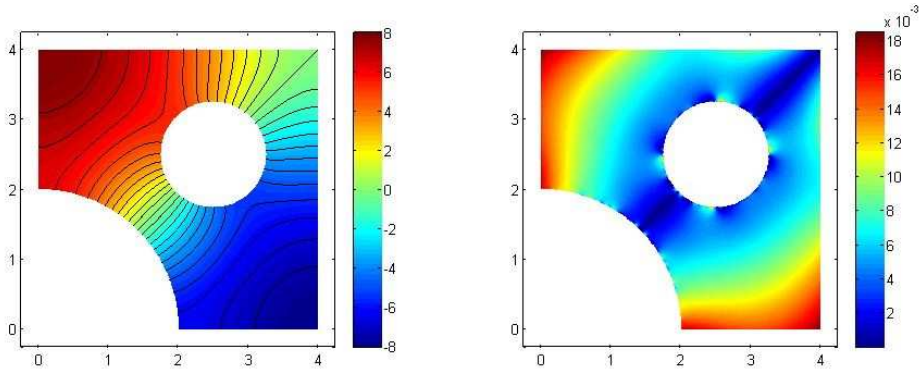,width=0.95\textwidth} 
 \caption{
Solution to Eq. (\ref{eqheat}) at time $t = 5$ in a more
complicated domain (left) and spatial distribution of the difference between the solutions obtained by the SSB method and using finite elements (right). Grid
resolution is $400 \times 400$ and $\xi = 0.05$.} \label{fig:heathole}
\end{center}
\end{figure}

An important advantage of the SSB method is that when the new formulation given by Eq. (\ref{equa4}) is used,
separate equations are not written for the boundaries, as the solution automatically adapts to satisfy the
boundary conditions on $\partial \Omega$, which results in a very simple computational implementation.
Alterations to the domain geometry therefore can be handled straightforwardly, without generating and
implementing additional boundary condition equations. As an example, we present in Fig. \ref{fig:heathole} the
solution to Eq. (\ref{eqheat}) in a more complicated geometry that combines both polar and Cartesian
coordinates. As we cannot obtain an analytic steady solution to the equation in this domain, we evolve it until time
$t = 5$ and compare our solution to one obtained using the \emph{finite elements} toolbox of {\sc Matlab}.
Since it is well known that spectral methods have higher accuracy than 
second-order finite elements, we obtained the finite-element 
solution using the finest possible grid that our machine could handle
(31152 nodes and 61440
triangles). At this resolution, as shown in Fig. \ref{fig:heathole}, the maximum difference between the solutions obtained by both
 methods is very small ($\approx 10^{-3}$). We expect that as the grid is refined for the finite elements method the solution will converge to the spectral solution, resulting in an even smaller difference error.

\subsection{The Allen-Cahn equation}

\begin{figure}
\begin{center}
\epsfig{file=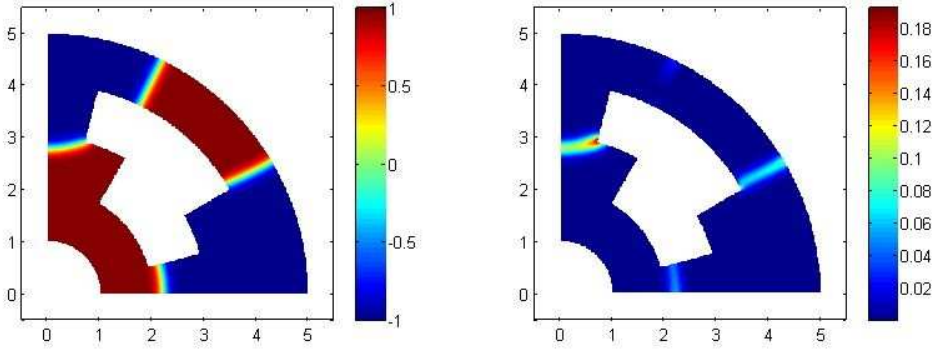,width=1.0\textwidth} 
\caption{
Solution of the Allen-Cahn equation (\ref{eq:allencahn}) at time $t = 65$ 
with $\epsilon = 10^{-2}$ using the SSB method (left) and errors when compared to the
solution obtained by a second-order finite-difference scheme in 
polar coordinates (right). Grid size is
$360 \times 360$ with $\xi = 0.05$ for the SSB method and 
grid resolution is $\Delta r=0.025$ and $\Delta \theta$=0.3 degrees for the polar finite-difference solution.}
\label{fig:allencahn}
\end{center}
\end{figure}

An important partial differential equation which arises in the modeling of the
formation and motion of phase boundaries is the
Allen-Cahn \cite{Allen} equation:
\begin{equation} \label{eq:allencahn}
\begin{array}{ccl}
\partial_t u = \epsilon^2 \Delta u - f(u), & & x \in \Omega \\
\partial_n u = 0, & & x \in \partial \Omega,
\end{array}
\end{equation}
where $\epsilon$ is a small positive constant and $f(u)$ is the
derivative of a potential function $W(u)$ that has two wells of
equal depth. For simplicity, we will assume that $W(u) =
(u^2-1)^2/4$, which makes $f(u) = u^3 - u$. In this manner, the
Allen-Cahn equation may be seen as a simple example of a nonlinear
reaction-diffusion equation. As explained in Ref.
\cite{Trefethen}, this equation has three fixed-point solutions,
$u=-1$, $u=0$ and $u = 1$. The middle state is unstable, but the
states $u = \pm 1$ are attracting, and solutions tend to exhibit
flat areas close to these values separated by interfaces that may
coalesce or vanish on long time scales, a phenomenon known as
metastability. Figure \ref{fig:allencahn} shows the solution of
the Allen-Cahn equation solved with Neumann boundary conditions on
an annulus with a $z$-shaped hole using the SSB method, with
$\epsilon = 0.01$. The annulus structure is given by $1 \le r \le
5$ and the z-hole is formed using radii at 2, 3, and 4 and angles
in steps of $15^\circ$ degrees (15, 30, 60 and $75^\circ$). For
initial conditions we have chosen two positive and two negative
Gaussian functions located in different parts of the domain:
\begin{equation} \label{eq:eqgaussian}
u=\sum_{i=1}^{n=4} (-1)^{n+1} \exp(-20 ((x-x_i)^2 + (y-y_i)^2)),
\end{equation}
where
\begin{equation}
\begin{tabular}{ll}
$x_1 = 1.5\cos(\pi/4)$,         &$y_1 = 1.5\sin(\pi/4)$, \\
$x_2 = 4\cos(\pi/12)$,          &$y_2 = 4\sin(\pi/12)$, \\
$x_3 = 4.5\cos(\pi/4)$,         &$y_3 = 4.5\sin(\pi/4)$, \\
$x_4 = 4\cos(11\pi/24)$,        &$y_4 = 4\sin(11\pi/24)$. \\
\end{tabular} \nonumber
\end{equation}

For comparison we also solved the equation using
a second order in space finite difference method in polar coordinates, since despite the
complexity of this geometry all the boundaries are parallel to the
axes in polar coordinates and so implementing no-flux boundary
conditions is straightforward. As
shown in Fig. \ref{fig:allencahn}, both the SSB and polar 
finite-difference solutions are in very good agreement, except in the
interface separating the two phases that lies just in one of the
corners. These larger but still overall small differences between both schemes are due to the sharp transition of the solution between -1 and 1 and the fact that the finite difference scheme approximate
spatial derivatives using lower-order accuracy than the spectral method.

\subsection{Reaction-diffusion equations and excitable media}

Models of excitable media form another significant class of
nonlinear parabolic partial differential equations and describe
systems as diverse as chemical reactions \cite{bz,bar},
aggregation of amoebae in the cellular slime mold Dictyostelium
\cite{Foerster}, 
%neural tissue \cite{Gorelova}, 
calcium waves
\cite{Marni}, and the electrical properties of neural \cite{HH}
and cardiac cells \cite{BR,LR}, among others. The equations of
excitable systems extend the Allen-Cahn equation by including one
or more additional variables that govern growth and decay of the waves.
Solutions of excitable media consist of excursions in state space
from a stable rest state and a return to rest, with the equations
describing the additional variables determining the time courses
of excitation and recovery. In spatially extended systems,
diffusive coupling allows excitation to propagate as nonlinear
waves, and in multiple dimensions complex patterns can be formed,
including two-dimensional spiral waves
\cite{bz,bar,Foerster,Marni,Davidenko} and their
three-dimensional analogs, scroll waves
\cite{gray,Fenton2}. Well-known examples of excitable
media equations include the Hodgkin-Huxley \cite{HH} model of
neural cells and its generalized simplification, the
FitzHugh-Nagumo \cite{FHN} model.

\begin{figure}
\begin{center}
\epsfig{file=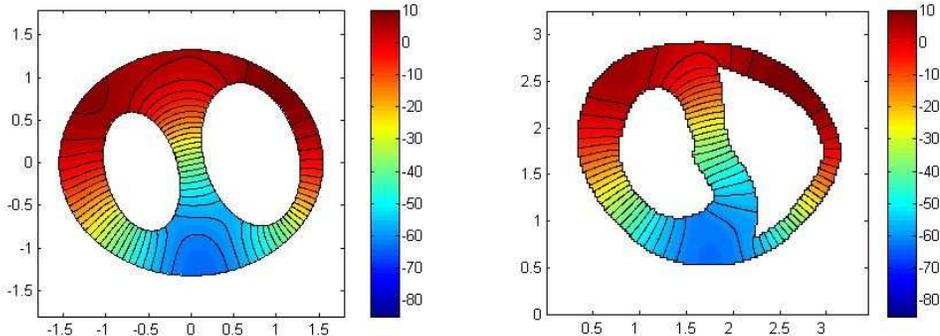,width=1.0\textwidth}
\caption{
Propagating wave of electrical potential in an idealized (left) and a realistic 
(right) slice of ventricular tissue using the SSB method. Grid size is 
$400 \times 400$ and $\xi = 0.025$ in both cases. 
The
color code denotes tissue voltage in mV, where 
red corresponds to cells with higher potential (depolarization) and blue to cells coming back to resting state (repolarization). Following Ref. \emph{\cite{Fenton2}}, parameters for these simulations in Eqs. \ref{eq:hearteqs} are $\tau_d = 0.25$, $\tau_r = 50$, $\tau_{\mathrm{si}} = 45$, $\tau_0 = 8.3$,  $\tau_v^{+} = 3.33$, $\tau_{v1}^{-} = 1000$, $\tau_{v2}^{-} = 19.2$, $\tau_w^{+} = 667$, $\tau_w^{-} = 11$, $u_c = 0.13$, $u_v = 0.055$ and $u_c^{\mathrm{si}} = 0.85$.
} \label{fig:heart}
\end{center}
\end{figure}

The dynamics of wave propagation in excitable media has been
studied extensively in regular domains. However, the complex geometry inherent
to some systems, such as the heart, often can have an essential influence on
 wave stability and dynamics \cite{Fenton3}. This fact, combined with the need for high-order accuracy to resolve the
sharp wave fronts characteristic of cardiac models, should make the SSB 
method a useful tool for studying electrical
waves in realistic heart geometries. Figure \ref{fig:heart} shows
an example of a propagating wave of action potential in both
an
idealized (left) and a realistic (right) slice \cite{Andrew} of ventricular tissue
using the SSB method and a phenomenological ionic cell model \cite{Fenton2,Fenton3}
with equations of the form
\begin{subequations}
\label{eq:hearteqs}
\begin{eqnarray}
& & \partial_t u({\bf x},t) = \nabla \cdot (\boldsymbol{D}\nabla
u)-J_{fi}(u,v)-J_{so}(u)-J_{si}(u,w) \\
& &\partial_t v({\bf x},t)  =  \Theta(u_c-u)(1-v)/\tau_v^{-}(u) - \Theta(u-u_c)v/\tau_v^{+} \\
& &\partial_t w({\bf x},t)  =  \Theta(u_c-u)(1-w)/\tau_w^{-} - \Theta(u-u_c)w/\tau_w^{+}  \\
& &J_{fi}(u,v)  =   -\frac{v}{\tau_d} \Theta(u-u_c) (1-u)(u-u_c) \\
& &J_{so}(u)  =   \frac{u}{\tau_0} \Theta(u_c-u) + \frac{1}{\tau_r} \Theta(u-u_c) \\
& &J_{si}(u,w)  =  -\frac{w}{2\tau_{\mathrm{si}}}(1+\tanh[k(u-u_c^{\mathrm{si}})]) \\
& &\tau_v^{-}(u)  =  \Theta(u-u_v)\tau_{v1}^{-} + \Theta(u_v-u)\tau_{v2}^{-}
\end{eqnarray}
\end{subequations}
 where u is the membrane potential; $J_{fi}, J_{so}$ and $J_{si}$ are
phenomenological currents; $v$ and $w$ are ionic gate variables; and $\boldsymbol{D}$ is the diffusion tensor (isotropic for these simulations, with value $D=1$ $\mathrm{cm}^2/\mathrm{s}$). In all these formulas, $\Theta(x)$ is the standard Heaviside step function defined by $\Theta(x)=1$ for $x \geq 0$ and $\Theta(x) = 0$ for $x < 0$, and the set of parameters of the model are chosen to reproduce different cellular dynamics measured experimentally.

\section{Discussion and future work}

In this paper we have presented a new method
for implementing homogeneous Neumann boundary conditions using
spectral methods for several problems of general
interest. The spectral smoothed boundary method offers several
advantages over finite-difference and finite-element
alternatives. Because ghost cells are not needed, the 
implementation of boundary conditions requires less coding than 
finite-difference
stencils, and in addition the spatial derivatives are represented
with higher accuracy. The use of simple Cartesian
grids makes the SSB method easier to use with multiple domain
shapes than finite elements, since \emph{grid generation
is not necessary}. Furthermore, the use of FFT routines in the
SSB method ensures efficiency and also makes extension
of the method to three dimensions straightforward using
well-established routines. Since the method is directly based on the FFT
it is very simple to implement on high performance computers by using 
 native parallel or vector FFT libraries.

The most significant limitation of the SSB method is that the 
error of the method depends directly on the ratio of the width of the smoothed boundary $\xi$ to the spatial step $\Delta x$, what implies that using uniform grids the number of points in the discretization needs to be large.
This limitation is perhaps not important for certain
classes of problems in which the solution contains steep wavefronts or other sharp features that require a fine
spatial resolution to correctly reproduce the dynamics of the system, such as electrical waves in cardiac tissue or shock waves in fluid mechanics. However, the adequate reduction of error in domains with irregular boundaries using this
method may require an increase in spatial resolution of a factor of 10 or 20 in each direction of the mesh for problems with smooth behavior like the heat equation with slowly varying sources compared to what is typically needed to obtain the same accuracy without complex boundaries.

Thus, an important future extension of this work is 
to improve the performance of the SSB method for problems
that do not track features with sharp spatial gradients.
For instance, if the boundary is stationary, it could be useful to use a non-uniform grid with extra resolution along
the boundaries combined with a non-uniform fast Fourier
transform (NFFT) to calculate the spatial derivatives.
However, if the boundary moves over time, it might be
more efficient to use a fine spatial discretization
than to keep track of the boundary for such problems. Other
planned future work includes properly handling complex
anisotropies in the diffusion matrices such as those found 
in cardiac muscle, examining
whether the method can be used to satisfy other types
of boundary conditions including Dirichlet and Robin, and
implementing non-stationary boundaries.

In conclusion, we have presented and analyzed a new numerical
method which imposes homogeneous Neumann boundary conditions
in complex geometries using spectral methods. We have
used this method to solve different partial differential
equations in domains with irregular boundaries and have found
 good agreement with the exact analytical solutions when such
solutions can be obtained. Along with the overall
advantage of allowing domains of different shapes to be considered
with spectral methods in a very simple way, this method also offers highly
accurate discretizations of spatial derivatives,
ease of implementation, 
straightforward extension to 3D, and
applicability to a wide variety of equations.
Moreover,  SSB codes need not change
to implement different geometries since all the information on
the geometry is contained in the function $\phi^{(\xi)}$, with the additional advantage that
 this function is
easy to generate and, unlike finite element methods, does not require the use of special
software for grid generation.

\hbox{}

\textbf{Acknowledgments.} We would like to
thank Elizabeth M. Cherry for useful discussions and valuable
comments on the manuscript. A. Bueno-Orovio also would like to thank the Physics Department at Hofstra University for their hospitality during his visit there (July - September, 2003). 
We also acknowledge the National Biomedical
Computation Resource (NIH Grant P41RR08605, USA). 

\bibliographystyle{siam}

\begin{thebibliography}{20}

\bibitem{Fornberg}
{\sc B. Fornberg}, {\em A practical guide to pseudospectral methods}, Cambridge University Press, Cambridge,
1996.

\bibitem{Trefethen}
{\sc Ll. N. Trefethen}, {\em Spectral methods in MATLAB}, SIAM, Philadelphia, 2000.

\bibitem{Sanzserna}
{\sc J.~M. Sanz-Serna}, {\em Fourier techniques in numerical methods for evolutionary problems}, in
 Proceedings of the Third Granada Seminar on Computational Physics, 1994,
 P.~L. Garrido and J. Marro, eds.,
Lecture Notes in Phys. 448, Springer-Verlag, Berlin, 1995, pp.~145--200.

\bibitem{Canuto}
{\sc C. Canuto, M. Y. Hussaini, A. Quarteroni and T. A. Zang}, {\em Spectral methods in fluid dynamics},
Springer-Verlag, Berlin, 1988.

\bibitem{Gottlieb}
{\sc  D. Gottlieb, S. A. Orszag}, {\em Numerical analysis of spectral methods: theory and applications},
SIAM, Philadelphia, 1977.

\bibitem{Roger}
{\sc R. Peyret}, {\em Spectral methods for incompressible viscous flow}, Springer-Verlag, New York, 2002.

\bibitem{solitones}
{\sc I. D. Iliev, E. Kh Khristov and K. P. Kirchev}, {\em Spectral methods in soliton equations}, CRC Press,
1994.

\bibitem{BunYuGuo}
{\sc Guo Ben-Yu}, {\em Spectral methods and their applications}, World Sientific, 1998.

\bibitem{Orszag}
{\sc S. A. Orszag}, {\em Spectral methods for problems in complex geometries}, J.
 Comput. Phys., \textbf{37}, 70 (1980).

\bibitem{Korczak}
{\sc K. Z. Korczak, A. T. Patera}, {\em An isoparametric spectral element method for solution of the
Navier-Stokes equations in complex geometry}, J. Comput. Phys., \textbf{62}, 361 (1986).

\bibitem{KarRap98}
{\sc A. Karma and W.-J. Rappel}, {\em  Quantitative phase-field modeling of dendritic growth in two and three
dimensions}, Phys. Rev. E, \textbf{57}, 4323 (1998).

\bibitem{Foletal99a}
{\sc R. Folch, J. Casademunt, A. Hern\'andez-Machado and L. Ram\'{\i}rez-Piscina}, {\em  Phase-field model for
Hele-Shaw flows with arbitrary viscosity contrast. I. Theoretical approach}, Phys. Rev. E, \textbf{60}, 1724
(1999).

\bibitem{Foletal99b}
{\sc R. Folch, J. Casademunt, A. Hern\'andez-Machado and L. Ram\'{\i}rez-Piscina}, {\em  Phase-field model for
Hele-Shaw flows with arbitrary viscosity contrast. II. Numerical study}, Phys. Rev. E, \textbf{60}, 1734 (1999).

\bibitem{Karetal01}
{\sc A. Karma, D. Kessler and H. Levine}, {\em Phase-field model of mode III dynamic fracture}, Phys. Rev.
Lett., \textbf{87}, 045501 (2001).

\bibitem{BibMis03}
{\sc T. Biben and C. Misbah}, {\em  Tumbling of vesicles under shear flow within an advected-field approach},
Phys. Rev. E, \textbf{67}, 031908 (2003).

\bibitem{Casademunt}
{\sc R. Gonzalez-Cinca, R. Folch, R. Benitez, L. Ramirez-Piscina, J. Casademunt, A. Hernandez-Machado}, {\em
Phase-field models in interfacial pattern formation out of equilibrium}, \texttt{arxiv.org/cond-mat/0305058} to
appear in Advances in Condensed Matter and Statistical Mechanics, ed. by E. Korucheva and R. Cuerno, Nova
Science Publishers.

\bibitem{FFTW}
{\sc M. Frigo and S.~G. Johnson}, {\em FFTW: An adaptive software architecture for the FFT}, in Proceedings of
the IEEE International Conference on Acoustics, Speech, and Signal Processing, IEEE, Seattle, WA, 1381 (1998).

\bibitem{Fenton}
{\sc F. H. Fenton, A. Karma and W.-J. Rappel}, {\em Modeling wave propagation in anatomical heart models using
the phase-field method}, (preprint).

\bibitem{Strang}
{\sc G. Strang}, {\em On the construction and comparison of difference schemes}, SIAM J. Num. Anal., \textbf{5},
506 (1968).

\bibitem{Allen}
{\sc S. M. Allen and J. W. Cahn}, {\em A microscopic theory for antiphase boundary motion and its application to
antiphase domain coarsening,}, Acta Metall. Mater., \textbf{27}, 1085 (1979).

\bibitem{bz}
{\sc A. N. Zaikin and A. M. Zhabotinsky}, {\em Concentration wave propagation in two-dimensional liquid-phase
self-oscillating system}, Nature, \textbf{225}, 535 (1970).

\bibitem{bar}
{\sc M. Bar, Ch. Zlicke, M. Eiswirth and G. Ertl}, {\em Theoretical modeling of spatiotemporal self-organization
in a surface catalyzed reaction exhibiting bistable kinetics}, J. Chem. Phys, \textbf{96}, 8595 (1992).

\bibitem{Foerster}
{\sc P. Foerster, S. C. Muller and B. Hess}, {\em Curvature and spiral geometry in aggregation patterns of
Dictyostelium discoideum}, Development, \textbf{109}, 11 (1990).

\bibitem{Marni}
{\sc M. E. Harris-White, S. A. Zanotti, S. A. Frautschy and A. C. Charles}, {\em Spiral intracellular calcium
waves in hippocampal slice cultures}, J. Neurobiol., \textbf{79}, 1045 (1998).

\bibitem{HH}
{\sc A. L. Hodgkin and A. F. Huxley}, {\em A quantitative description of membrane current and its application to
conduction and excitation in nerve}, J. Physiol., 117, 500, (1952).

\bibitem{BR}
{\sc G. W. Beeler and H. Reuter.}, {\em Reconstruction of the action potential of ventricular myocardial
fibres}, J. Physiol., \textbf{268}, 177 (1977).

\bibitem{LR}
{\sc C. Luo and Y. Rudy}, {\em A model of the ventricular cardiac action potential. Depolarization,
repolarization, and their interaction}, Circ. Res., \textbf{68}, 1501 (1991).

\bibitem{Davidenko}
{\sc  J. M. Davidenko, A. M. Pertsov, R. Salomonsz, W. Baxter and J. Jalife}, {\em Stationary and drifting
spiral waves of excitation in isolated cardiac muscle}, Nature, \textbf{355}, 349 (1992).

\bibitem{gray}
{\sc R. A. Gray, A. M. Pertsov and J. Jalife}, {\em Spatial and temporal organization during cardiac
fibrillation}, Nature, \textbf{392}, 75 (1998).

\bibitem{Fenton2}
{\sc F. Fenton and A. Karma}, {\em Instability of electrical vortex filament and wave turbulence in thick
cardiac muscle}, Phys. Rev. Lett., \textbf{81}, 481 (1998).

\bibitem{FHN}
{\sc R. FitzHugh}, {\em Impulses and physiological states in theoretical
 models of nerve membranes}, Biophys. J., \textbf{1}, 445 (1961).

\bibitem{Fenton3}
{\sc F. H. Fenton, E. M. Cherry, H. M. Hastings and S. J. Evans}, {\em Multiple mechanisms of spiral wave
breakup in a model of cardiac electrical activity}, Chaos, \textbf{12}, 852 (2002).

\bibitem{Andrew}
{\sc F. J. Vetter, A. D. McCulloch}, {\em Three-dimensional analysis of regional cardiac function: a model of
rabbit ventricular anatomy}, Prog. Biophys. Mol. Biol. \textbf{69}, 157 (1998).

\end{thebibliography}

\end{document}